\documentclass{amsart}

\usepackage[all,hyperref]{preprint_style}
\usepackage{amsmath,amssymb,amsbsy,amsthm,amsfonts}
\usepackage[foot]{amsaddr}

\usepackage[T1]{fontenc}

\usepackage{caption}
\usepackage[caption=false]{subfig}

\usepackage{mathtools}
\usepackage{dsfont}
\usepackage{a4wide}
\usepackage{bm}

\usepackage{graphicx}
	\DeclareGraphicsExtensions{.pdf,.png,.jpg}
\usepackage{pgfplots}
	\pgfplotsset{compat=1.3}
\usepackage{placeins}
\usepackage{tabularx}
\usepackage{enumitem}
\usepackage{tikz}

\numberwithin{equation}{section}
 \makeatother
\def\ll{\left\langle}

\def\fp{\,{:}\,}

\def\tens#1{\pmb{\mathsf{#1}}}
\def\vec#1{\boldsymbol{#1}}

\def\sym{\mathop{\mathrm{sym}}\nolimits}

\def\loc{\mathop{\mathrm{loc}}\nolimits}

\def\tr{\mathop{\mathrm{tr}}\nolimits}

\def\diver{\mathop{\mathrm{div}}\nolimits} 

\def\bf{\vec{f}}

\def\bn{\vec{n}}

\def\bu{\vec{u}}
\def\bv{\vec{v}}
\def\bw{\vec{w}}

\def\bz{\vec{z}}

\def\BA{\tens{A}}

\def\BD{\tens{D}}

\def\BG{\tens{G}}

\def\BI{\tens{I}}

\def\BS{\tens{S}}
\def\BT{\tens{T}}

\def\Bsigma{\tens{\sigma}}

\def\tria{\mathcal{T}_h}
\def\faces{\Gamma_h}

\def\leftjump{\left[\!\left[}
\def\rightjump{\right]\!\right]}
\def\leftavg{\left\{\!\!\left\{}
\def\rightavg{\right\}\!\!\right\}}
\def\avg#1{\{\!\!\{ #1 \}\!\!\}}

\def\Re{\mathrm{Re}}
\def\Eu{\mathrm{Eu}}
\def\Ga{\mathrm{Ga}}

\def\Lsymtr{L^2_{\sym,\tr}(\Omega)^{d\times d}}
\def\Hone{H_0^1(\Omega)^d}
\def\Hdiv{H(\diver;\Omega)}
\def\Lmean{L^2_0(\Omega)}

\begin{document}

\title[Development of boundary layers in activated Euler fluids ]{Development of boundary layers in Euler fluids that on ``activation'' respond like Navier-Stokes fluids}

\thanks{J.~M\'{a}lek acknowledges the support of the project No. 20-11027X financed by the Czech Science Foundation (GA \v{C}R). J. M\'{a}lek is a member of the Ne\v{c}as Center for Mathematical Modelling. K~R.~Rajagopal thanks the Office of Naval Research for its support of this work.}

\author[Gazca-Orozco]{P.~A.~Gazca-Orozco$^\star$}
\address[$\star$]{Department of Applied Mathematics, University of Freiburg, Ernst–Zermelo–Straße 1, 79104 Freiburg, Germany}
\email{alexei.gazca@mathematik.uni-freiburg.de}

\author[Málek]{J.~Málek$^\dagger$}
\address[{$\dagger$}]{Faculty of Mathematics and Physics, Charles University, Sokolovská 83, 18675 Prague 8, Czech Republic}
\email{malek@karlin.mff.cuni.cz}

\author[Rajagopal]{K.~R.~Rajagopal$^\ddagger$}
\address[{$\ddagger$}]{Department of Mechanical Engineering,  
Texas A\&M University, College Station, TX 77845 USA}
\email{krajagopal@tamu.edu}

\date{\today}

\maketitle

\begin{abstract}
	\vspace{-6ex}
	We consider the flow of a fluid whose response characteristics change due the value of the norm of the symmetric part of the velocity gradient, behaving as an Euler fluid below a critical value and as a Navier-Stokes fluid at and above the critical value, the norm being determined by the external stimuli. We show that such a fluid, while flowing past a bluff body, develops boundary layers which are practically identical to those that one encounters within the context of the classical boundary layer theory propounded by Prandtl. Unlike the classical boundary layer theory that arises as an approximation within the context of the Navier-Stokes theory, here the development of boundary layers is due to a change in the response characteristics of the constitutive relation. We study the flow of such a fluid past an airfoil and compare the same against the solution of the Navier-Stokes equations. We find that the results are in excellent agreement with regard to the velocity and vorticity fields for the two cases.
\end{abstract}


\section{Introduction}

The main motivation behind the development of boundary layer theory is Prandtl’s cognizance \cite{PR.1904} that the effect of viscosity is restricted to a narrow region adjacent to a solid boundary past which fluids like air and water flow. Outside this narrow layer referred to as the “boundary layer” the fluid flows like an inviscid fluid, while even within the layer further approximations are made to simplify the Navier-Stokes equations, the equations being referred to as the “boundary layer equations”. A detailed exposition of boundary layer theory and an extensive bibliography related to the theory can be found in the book by Schlichting \cite{Schli}. An interesting issue in the development of boundary layers is the boundary layer thickness which delineates where the boundary layer equations ought to be enforced. Usually this is determined by matching the solutions for the boundary layer equations in a domain adjacent to the boundary and solving the Euler equations outside the domain and matching the two solutions to arrive at the boundary of the boundary layer. 

Within the context of the Navier-Stokes theory, boundary layers are a consequence of the inertial term, and the higher the Reynolds number, the more pronounced the boundary layer. This has led to the misconception that boundary layers manifest themselves only at sufficiently high Reynolds numbers. Such an understanding is incorrect. Boundary layers, in the sense that vorticity is confined to a specific region, can manifest themselves in nonlinear fluids even in the limit of zero Reynolds number (see Rajagopal \cite{KRR1995}, Mansutti and Rajagopal \cite{ManKRR}), and even at zero Reynolds number multiple deck structured boundary layers are possible with the effects of different physical quantities confined to such regions (see Rajagopal, Gupta and Wineman \cite{RajGW}). Put simply, boundary layers can manifest themselves due to an appropriate structure to the nonlinearities that arise in the governing equation. Of course, the boundary layers that arise in these different nonlinear fluids has little to do with the boundary layer that manifests itself due to the flow of a Navier-Stokes fluid past a bluff body. 

Implicit constitutive theories present a natural way to showcase the interesting circumstance wherein a bifurcation of the constitutive relations occurs, that is, there is a possibility that the body under consideration can be described by more than one constitutive relation based on some criterion for the selection of the appropriate constitutive relation (see Rajagopal and Wineman \cite{RW.1980} for the bifurcation of response in inelastic bodies, Rajagopal and Srinivasa \cite{RS.2004ii} for an explanation of solid to solid phase transition such as from Martensite to Austenite using the criterion of the maximality of the rate of entropy production, and Cichra and Průša \cite{Cichra2020} for a more recent development within a complete thermodynamic setting).


We are interested in a fluid that responds as an Euler fluid if the norm of the symmetric part of the velocity gradient is below a certain threshold, and behaves as a Navier-Stokes fluid when the norm of the symmetric part of the velocity gradient is at or above the threshold. That is, it is in the nature of the fluid to change its response characteristics based on the stimuli, which in this case is the extent of the shear-rate as a consequence of the conditions that the fluid is subject to at the boundaries of the flow domain or for that matter body forces acting on the fluid. We refer to such a fluid as an activated Euler fluid, a fluid constitutive relation considered by Blechta et al. ~\cite{BMR.2020}. Such a bifurcation of response characteristics is different from the bifurcation of solutions to the equations of motion that might correspond to the governing equations for a specific constitutive relation. These two kinds of bifurcation, that of the constitutive relation, and that for the solution to the equations of motion, are completely different in character and together offer a far richer class of solutions to describe a specific phenomenon whose explanation is being sought.

From the viewpoint of the mathematical properties of the solution of the governing partial differential equations, the solutions to the governing equations for the activated Euler fluid present much better characteristics than those presented by the Euler fluid in that one can establish global-in-time existence of weak solutions for large data for a variety of boundary conditions (see Blechta et al. \cite{BMR.2020} for details of the same).

The boundary layers that are produced adjacent to a solid surface and bluff bodies due to the flow of an activated Euler fluid are very similar in structure to the boundary layers that are created by the flow of a Navier-Stokes fluid. That is, the vorticity is confined to a narrow region adjacent to a solid wall, and outside of this boundary layer the solution is governed by that for an Euler fluid.  

The important difference between the classical boundary layer approach and using such an activated Euler fluid is that there is no necessity to match the solutions of two distinct sets of equations to determine the domain of application of the boundary layer equations. Moreover, there are no approximations as in the obtainment of the boundary layer equations; instead the full equations that stem from the balance of linear momentum for the constitutive relation are solved, and automatically the constitutive relation changes so that the equations reduce to the Euler equation for a very large flow domain with a more complicated equation to be solved only in a small region, as in the case of boundary layer theory.

It would be appropriate to mention that the governing equations that arise present challenging questions from the perspective of numerical analysis. The constitutive equation for the activated Euler fluid is explicit (the stress is a continuous function of the velocity gradient), and it may seem natural to apply any Navier-Stokes code as computational solver. However, as the governing equations are non-smooth and the nonlinearity generates a non-strictly monotone operator (see the equation \eqref{josef2} or \eqref{eq:actEuler_dichotomy}) we proceed first to regularise the problem. We exploit the dual relationship between our model and the Bingham model, both belonging to the class of non-smooth constitutive relations, in carrying out the regularisation. As a consequence, the regularised problem belongs to the class of incompressible implicitly constituted fluids (for the velocity $\bv$, the pressure $p$ and the stress $\BS$, see \eqref{eq:non_dimensional}). Consequently, instead of the otherwise natural $(\bv,p)$ formulation, a $(\bv,p,\BS)$ formulation becomes more natural.


In order to illustrate the efficacy of the constitutive relation in capturing a boundary layer similar in structure to that produced by the flow of a Navier-Stokes fluid, we solve the model problem of flow past an airfoil and compare the result for a Navier-Stokes fluid with that for an activated Euler fluid. We find that the solutions in the two cases are nearly the same.

The idea to split the domain occupied by a Navier-Stokes fluid into the region where the behavior of the fluid is sufficiently well described by the Euler fluid and the region where the viscous effects are  taking place has been exploited in designing numerical methods before. Here, in particular we refer to \cite{BCR.89,AP.93,GP.2017B}, see also our concluding remarks for further comments.

The organization of the paper is as follows. In the next section we introduce the constitutive relation for an activated Euler fluid. Viewing it as a dual constitutive equation to the one for a Bingham fluid model, within the framework of implicitly constituted incompressible fluid models, we introduce a regularised approximation. Section 3 is devoted to the formulation of the problem, its discretization and the development of the computational scheme. In section 4 we study computationally the flow of an activated Euler fluid past a bluff body and compare the results against those for a Navier-Stokes fluid. In the final section we provide some concluding remarks.

\section{Activated Euler fluid}

The governing equations for flows of implicitly constituted incompressible fluids (see Rajagopal \cite{Raj.2003, Raj.2006}, M\'{a}lek and Rajagopal \cite{MR.2005, MR.2007}) consist of the balance equations 
\begin{equation}
    \diver \bv = 0 \quad \textrm{ and } \quad \rho_* \frac{\textrm{d}\bv}{\textrm{d}t} = - \nabla p + \diver \BS + \rho_* \bm{f}\,, \label{eq:p1}
\end{equation}
and the constitutive equation
\begin{equation}
    \BG(\BS, \BD) = 0. \label{eq:p2}
\end{equation}
Here, $\bv$ is the velocity, $p$ the pressure, $\rho_*$ is a given constant density, $\bm{f}$ represents the density of the given external body borces, $\BS := \BT - \frac{1}{d}\tr(\BT)\BI$ is the deviatoric (traceless) part of the Cauchy stress $\BT$. 
Furthermore, $\BD:= \BD\bv = \tfrac{1}{2}(\nabla \bv + \nabla \bv^\top)$ is the symmetric part of the velocity gradient and $\BG$  stands for a given symmetric tensor-valued functions of two tensor variables. The equations \eqref{eq:p1} and \eqref{eq:p2} are considered in a $d$-dimensional flow domain $\Omega$, $d\ge 2$.
  
Recently, Blechta et al. \cite{BMR.2020} developed a systematic approach to study a class of implicit constitutive equations of the form \eqref{eq:p2}. As an outcome, a new class of models have been identified, namely fluids that behave as an Euler fluid until the (magnitude of the) shear rate does not exceed a certain critical value; once this happens the fluid responds as a Navier-Stokes fluids (or a more complex non-Newtonian fluid). We call fluids described in this way \emph{activated Euler fluids}. We briefly recall their ``derivation''.

We start by noticing that if $\BS$ and $\BD$ are related linearly in \eqref{eq:p2}, we obtain the constitutive equation for a Navier-Stokes fluid: 
\begin{equation*}
     \BS = 2\nu_* \BD,         
\end{equation*}
where $\nu_\star>0$ is the viscosity.

    
    Another very popular model belonging to the class \eqref{eq:p2} with broad applications is the Bingham model, which is usually described by the relations
    \begin{equation}\label{eq:Bingham_dichotomy}
  \renewcommand{\arraystretch}{1.2}
\left\{
	\begin{array}{ccl}
		|\BS|\leq \sigma_\star & \Longleftrightarrow & \BD = \bm{0}, \\
|\BS|> \sigma_* & \Longleftrightarrow & \BS = 2\nu_\star \BD + \sigma_\star\displaystyle\frac{\BD}{|\BD|},
\end{array}
\right.
\end{equation}
where $\sigma_\star>0$ is the activation (yield) stress. The Bingham fluid model exhibits a response that is \emph{non-linear}, \emph{monotone} but not strictly monotone, and continuous but \emph{non-smooth}. Moreover, viewing \eqref{eq:Bingham_dichotomy} as an $\BS$ vs. $\BD$ relation, this represents a \emph{multi-valued} mapping as for $\BD=\bm{0}$ there is infinite many admissible values of $\BS$ satisfying \eqref{eq:Bingham_dichotomy}. Changing however the viewpoint, and looking at \eqref{eq:Bingham_dichotomy} as a $\BD$ vs. $\BS$ relation we observe that $\BD$ is a \emph{single-valued} function of $\BS$ that has, for instance, the following explicit form: 
\begin{equation}
    \label{josef1}
    \BD = \frac{1}{2\nu_\star}[\, |\BS| - \sigma_\star \,]_+\frac{\BS}{|\BS|}\,,
\end{equation}
where $[z]_+$ denotes a positive part of $z$, i.e. $[z]_+:= \max\{z,0\}$. Interchanging the role of $\BS$ and $\BD$ in \eqref{josef1} and setting $\alpha_\star:= \frac{1}{2\nu_\star}$ (fluidity) and relabelling $\tau_\star:=\sigma_\star$ we obtain 
\begin{equation}
    \label{josef2}
    \BS = \alpha_\star [\, |\BD| - \tau_\star\, ]_+\frac{\BD}{|\BD|}.
\end{equation}
This is the aforementioned \emph{activated Euler fluid model} that we investigate in this study. Note that \eqref{josef2} can 
alternatively be described in terms the dichotomy:
\begin{equation}\label{eq:actEuler_dichotomy}
  \renewcommand{\arraystretch}{1.2}
\left\{
	\begin{array}{ccc}
		|\BD|\leq \tau_\star & \Longleftrightarrow & \BS = \bm{0}, \\
|\BD|> \tau_* & \Longleftrightarrow & \BD = \alpha_\star \BS + \displaystyle\frac{\tau_\star}{|\BS|}\BS .
	\end{array}
\right. 
\end{equation}

The constitutive relation \eqref{josef2} (similar to  \eqref{josef1}) is nonlinear, non-smooth (the relations are non-differentiable at the activation point) and monotone but not strictly monotone. In developing a suitable numerical method we benefit from the mutually dual relationship between the Bingham model \eqref{josef1} on one side, and activated Euler model \eqref{josef2} on the other side. 	More precisely, we follow in our approach regularisation strategies analogous to those popular in the analysis and computation of Bingham fluid flow. 
  
	Interestingly, for the Bingham fluid model the $\varepsilon$-regularization leads to a model where the stress $\BS$ is the function of $\BD$, such as for example the relation\footnote{This type of regularization is sometimes called the Bercovier-Engelman regularisation of the Bingham constitutive relation, see \cite{BE.1980}}
\begin{equation}\label{eq:Bingham_regularised}
\BS 
= \nu_g(|\BD|)\BD := 2\nu_\star \BD
+ \frac{\sigma_\star}{\sqrt{|\BD|^2 + \varepsilon_\star^2}}\BD,
\end{equation}
where $\nu_g$ is a generalised viscosity; note that in this case one can either employ a $(\bv,p)$ or $(\bv,p,\BS)$ formulation. 

In contrast, an anagolous $\varepsilon_\star$-regularization for the activated Euler model takes the form  
\begin{equation}\label{eq:actEuler_regularised}
\BD 
= \alpha_g(|\BS|)\BS :=
 \alpha_\star \BS
+
\frac{\tau_\star}{\sqrt{|\BS|^2 + \varepsilon_\star^2}}\BS,
\end{equation}
where $\alpha_g$ now represents a generalised fluidity; 
note that this expression does not allow one to insert $\BS$
in the balance of linear momentum and one is forced to work with a $(\bv,p,\BS)$-formulation. 

\section{Description of the problem - regularization and weak formulation}
We are interested in solving the following system of PDEs:
\begin{equation}
\begin{aligned}\label{eq:PDE}
			\rho_\star \diver (\bv \otimes\bv) 
			- \diver\, \BS \,+\, \nabla p 
			&= \rho_\star \bm{f}_\star
			\quad & &\text{ in }\Omega,\\
			\diver\bv &= 0\quad & &\text{ in }\Omega,\\
            \BS &= \alpha_* \left[\, |\BD| - \tau_\star \,\right]_{+} \frac{\BD}{|\BD|} \quad & &\text{ in }\Omega, \\
			\bv &= \bm{0} && \text{ on }\partial\Omega,
	\end{aligned}
\end{equation}
where $\alpha_\star>0$ is the fluidity and $\tau_\star\geq 0$ is the activation parameter.

In order to make the application of Newton's method feasible, we regularise the constitutive relation through the expression \eqref{eq:actEuler_regularised}. This means that the problem that is solved in this study is the following: for given positive parameters $\rho_*$, $\alpha_*$, $\tau_*$ and for a given function $\bm{f}_\star:\Omega \to \mathbb{R}^d$ we look for $(\bv,p,\BS):\Omega \to \mathbb{R}^d \times\mathbb{R}\times \mathbb{R}^{d\times d}$ satisfying 
\begin{equation}
\begin{aligned}\label{eq:PDEeps}
			\rho_\star \diver (\bv \otimes\bv) 
			- \diver\, \BS \,+\, \nabla p 
			&= \rho_\star \bm{f}_\star
			\quad & &\text{ in }\Omega,\\
			\diver\bv &= 0\quad & &\text{ in }\Omega,\\
            \BD &= \alpha_\star \BS + \frac{\tau_\star}{\sqrt{|\BS|^2 + \varepsilon_\star^2}}\BS \quad & &\text{ in }\Omega, \\
			\bv &= \bm{0} && \text{ on }\partial\Omega,
	\end{aligned}
\end{equation}

To find a non-dimensional version of the system, note that a characteristic shear-rate $\gamma_c = \frac{U_c}{L_c}$ (here $U_c$ and $L_c$ are the characteristic velocity and length scales, respectively) determines uniquely a characteristic stress $\sigma_c$ via the relation:
\begin{equation}
\gamma_c = \alpha_g(\sigma_c)\sigma_c.
\end{equation}
This allows us to define the characteristic fluidity as $\alpha_c := \alpha_g(\sigma_c)\sigma_c$. Following the traditional non-dimensionalisation procedure we arrive at the system:
\begin{equation}\label{eq:non_dimensional}
\begin{gathered}
\diver \bv = 0, \\
\mathrm{Re}(\diver(\bv\otimes \bv)) - \diver\BS + \nabla p = \mathrm{Ga} \bm{f}, \\
\BD = \alpha \BS + \frac{\mathrm{Eu}}{\sqrt{|\BS|^2 + {\varepsilon}^2}} \BS, \\
\mathrm{Re} = \rho_\star \alpha_c L U,
\quad
\alpha := \frac{\alpha_\star}{\alpha_c},
\quad
\mathrm{Eu} := \frac{\tau_\star L_c}{U_c},
\quad
\varepsilon^* := \frac{\varepsilon_\star \alpha_c L_c}{U_c},
\quad
\Ga := \frac{\alpha_c \rho_\star |\bm{f}_\star|L_c^2}{U_c}.
\end{gathered}
\end{equation}

	Before introducing a weak formulation of the system \eqref{eq:non_dimensional}, we define some useful function spaces; by $L^2(\Omega)$ we will denote the space of square-integrable functions on $\Omega$. Moreover, we define:
	\begin{gather*}
		H^1_0(\Omega) := \{\bw\in L^2(\Omega)^d \, :\, \nabla\bw\in L^2(\Omega)^{d\times d}\text{ and }\bw|_{\partial\Omega}=0 \}, \\
		L^2_{\sym,\tr}(\Omega)^{d\times d} := \{\BA \in L^2(\Omega)^{d\times d} \, :\, \tr(\BA) = 0\text{ and } \BA^\top = \BA\}, \\
		L^2_0(\Omega) := \left\{q\in L^2(\Omega) \, :\, \int_\Omega q = 0\right\}, \\
		H(\diver;\Omega) := \{\bw \in L^2(\Omega)^d \, :\, \diver\bw \in L^2(\Omega)\}.
	\end{gather*}

In a weak formulation of the system \eqref{eq:non_dimensional}, we then look for a triplet $(\BS,\bv,p)\in \Lsymtr\times \Hone \times \Lmean$ such that:
 \begin{equation}\label{eq:weak_form}
\begin{aligned}
													\int_\Omega \BD(\bv)\fp \BA 
													&=
													\int_\Omega \left(\alpha+ \frac{\Eu}{\sqrt{|\BS|^2 + \varepsilon^2}} \right) \BS\fp \BA
													&&\forall\, \BA\in \Lsymtr.\\
\int_\Omega \BS \fp \BD(\bw)
- \Re \int_\Omega (\bv &\otimes \bv)\fp \BD(\bw)
- \int_\Omega p\diver  \bw 
=
\Ga \int_\Omega \bm{f}\cdot \bw
											 &&
 \forall \,\bw\in \Hone. \\
											 &-\int_\Omega \diver \bv\, q = 0
													&& 
													\forall\, q\in L^2(\Omega)
\end{aligned}
\end{equation}
This system has a solution and we know for instance that for a given velocity $\bv$, the stress $\BS$ and pressure $p$ are uniquely defined (and no other velocities are associated to $\BS$); for $\varepsilon=0$ the model still has a solution but some uniqueness properties may get lost; see \cite{BMR.2020} for more details.
In fact a rigorous mathematical foundation (in terms of long-time and large-data existence of a weak solution and its properties) has been established in \cite{BMR.2020} for steady and unsteady (internal) flows, including no-slip and Navier's boundary conditions; this is e.g.\ in contrast to the ad-hoc activation models used in \cite{BCR.89,AP.93}.

\section{Discretisation and Solvers}
Let $\{\mathcal{T}_h\}_h$ be a sequence of shape-regular simplicial triangulations of $\Omega$. We will make use of the following finite element spaces for the discretisation of the stresses, velocities and pressures, respectively:
\begin{gather*}
\Sigma_h := \{\Bsigma_h \in \Lsymtr \,\mid \, \Bsigma_h|_K \in \mathbb{P}_0(\tria)^{d\times d} \text{ for all }K\in \tria \},\\
V_h := \{\bv \in \Hdiv \,\mid \, \bv|_K \in \mathbb{BDM}_1(\tria) \text{ for all }K\in \tria \},\\
Q_h := \{q \in \Lmean \,\mid \, q|_K \in \mathbb{P}_0(\tria) \text{ for all }K\in \tria \}.
\end{gather*}
Here $\mathbb{P}_0(\tria)$ and $\mathbb{BDM}_1(\tria)$ denote the space of piecewise constant functions and the Brezzi--Douglas--Marini finite element space of lowest order, respectively.
The choice of piecewise constant stresses is better suited to models containing activation parameters; for higher order elements the constitutive model cannot be evaluated exactly (i.e.\ $\alpha_g(|\BS_h|)\BS_h \not\in \Sigma^h$), which can lead to numerical issues (in a sense, $\alpha_g^{-1}$ approximates discontinuous behaviour, so higher-order elements can lead to oscillations); see \cite{TRFW.2018} for similar considerations in the case of the flow of a Bingham fluid.
This can be alleviated with adaptive mesh refinement, but since this is outside the scope of this work, we stay with the lower order discretisations stated above. 

The velocity space based on Brezzi-Douglas-Marini elements has the important property that discretely divergence-free functions are in fact pointwise divergence-free \cite{BBF.2013}; i.e.\ if $\bw_h \in V^h$, we have
\begin{equation*}
\int_\Omega \diver\bw_h \, q_h = 0 \quad
\forall\, q_h\in Q^h
\qquad
\Longrightarrow
\qquad
\diver \bw_h(x) = 0
\text{ for all }x\in \Omega.
\end{equation*}
This means that this discretisation preserves more faithfully the properties of the continuous system, and in particular leads to pressure-robustness of the scheme, a property whose importance has been recognised in recent years \cite{JLMNR.2017}. A drawback of this choice is that the space is no longer conforming, i.e.\ $V^h\not\subset \Hone$, and it is necessary to penalise the jumps across facets to recover a function without jumps in the limit $h\to 0$. In particular, the space is endowed with the norm
\begin{equation}\label{eq:velocity_norm}
\|\bw_h\|_{1,h} :=
\left(
\|\nabla_h \bw_h\|^2_{L^2(\Omega)}
+
\|h_\Gamma^{-1} \leftjump \bw_h\otimes \bn \rightjump\|^2_{L^2(\faces)}
\right)^{1/2},
\end{equation}
where $\nabla_h$ denotes the broken gradient ($(\nabla_h w_h)|_K = \nabla w_h|_K$), $h_\Gamma$ is the local face-size function defined on the mesh facets ($h_\Gamma \colon \faces \to R$, $h_\Gamma|_F = h_F$), and $\leftjump \bw_h\otimes \bn \rightjump := \bw_h^+\otimes \bn^+ + \bw_h^- \otimes \bn^-$ denotes the jump across a facet with normal vector $\bn$. As usual in DG notation, we will denote the average of a function across a facet as $\leftavg \phi\rightavg := \tfrac{1}{2}(\phi^+ + \phi^-)$, and set $\leftjump \bw_h \otimes \bn \rightjump = \bw_h \otimes \bn$ and $\leftavg \phi \rightavg = \phi$ on facets belonging to the boundary of $\Omega$.

The velocity space $V_h$ pairs with the stress and pressure spaces in a stable manner, meaning that the inf-sup conditions that ensure the well-posedness of the discrete system are satisfied (see e.g.~\cite{FGS.2020}); more precisely, one has 
\begin{subequations}
	\begin{alignat}{2}
\gamma_1
\|\bw_h\|_{1,h}
\leq
\sup_{\Bsigma_h\in \Sigma^h}
\frac{\int_\Omega \Bsigma_h\fp \BD\bw_h }{\|\Bsigma_h\|_{L^2(\Omega)}}
&+
\|h_\Gamma^{-1} \leftjump \bw_h\otimes \bn \rightjump\|_{L^2(\faces)}
\qquad
&\forall\, \bw_h\in V^h, \label{eq:infsup_Tv}
\\
\gamma_2 \|q_h\|_{L^2(\Omega)}
\leq
\sup_{\bw_h \in V^h}
&\frac{\int_\Omega q_h \diver \bw_h}{\|\bw_h\|_{1,h}}
&\forall\, q_h\in Q^h,  \label{eq:infsup_vp}
\end{alignat}
\end{subequations}
for two positive constants $\gamma_1,\gamma_2$, that are independent of the mesh size;
regarding the first condition \eqref{eq:infsup_Tv}, noting that $\BD(V_h) \subset \Sigma_h$ one has for any $\bw_h\in V_h$:
\begin{equation*}
\sup_{\Bsigma_h\in \Sigma^h}
\frac{\int_\Omega \Bsigma_h\fp \BD\bw_h }{\|\Bsigma_h\|_{L^2(\Omega)}}
=
\|\BD \bw_h\|_{L^2(\Omega)}.
\end{equation*}
This means that the condition \eqref{eq:infsup_Tv} is nothing but a reformulation of the Korn/Poincaré inequality for DG spaces, whose proof can be found e.g.\ in \cite{BdPG.2019}. 
The second condition \eqref{eq:infsup_vp} is classical in the analysis of isochoric flow, see e.g.~\cite{BBF.2013}.

Employing fluxes associated with the Local Discontinuous Galerkin method, in the discrete formulation we look for a triplet $(\BS_h,\bv_h,p_h)\in \Sigma^h \times V^h \times Q^h$ such that:
 \begin{equation}\label{eq:discrete_form}
\begin{gathered}
													\int_\Omega \BD_h(\bv_h)\fp \BT_h 
													-
													\int_{\faces} \leftjump \bv_h \otimes \bn \rightjump \fp \avg{\BT_h}
													-
													\int_\Omega \alpha_g(|\BS_h|)\BS_h \fp \BT_h = 0
													\hspace{5ex}\forall\, \BT_h\in \Sigma^h.\\
													\hspace{-12ex}
\int_\Omega \BS_h \fp \BD_h(\bv_h)
-
\int_{\faces} \avg{\BS_h} \fp \leftjump \bv_h \otimes \bn \rightjump
													+ \delta\int_{\faces} h_\Gamma^{-1} \leftjump \bv_h \otimes \bn \rightjump \fp \leftjump \bv_h\otimes \bn \rightjump  \hspace{4ex} \\
													\hspace{12ex}
- \Re \int_\Omega (\bv_h \otimes \bv_h)\fp \BD_h(\bv_h)
- \int_\Omega p_h\diver  \bv_h
=
\Ga \int_\Omega \bm{f}\cdot \bv_h
\hspace{4ex}
											 \forall \,\bv_h\in V^h. \\
											 -\int_\Omega \diver \bv_h\, q_h = 0
											 \hspace{8ex}
													\forall\, q_h\in Q^h.
\end{gathered}
\end{equation}

In order to solve the nonlinear discrete problem \eqref{eq:discrete_form}, Newton's method will be employed, and in each iteration we employ a multigrid method based on augmented Lagrangian preconditioners inspired by \cite{FG.2020,LFM.2022}.
	To be more precise, after linearisation, the discrete system \eqref{eq:discrete_form} can be written in the following form at each Newton iteration:
	\begin{equation}\label{eq:linearisation}
		\begin{bmatrix}
			A & B^\top \\
			B & 0
		\end{bmatrix}
		\begin{bmatrix}
\bz_h \\ p_h
		\end{bmatrix}
		=
		\begin{bmatrix}
			\bm{F} \\
			g
		\end{bmatrix},
	\end{equation}
	where $\bz_h := (\BS_h, \bu_h)^\top$, $A$ is the (linearised) stress-velocity block, and $B$ is the divergence operator acting on $V^h$; $\bm{F}$ and $g$ then represent the appropriate right-hand-sides arising from the Newton linearisation. The main strategy is to employ a block preconditioner based on the approximation:
	\begin{equation}\label{eq:block_factors}
		\begin{bmatrix}
			A & B^\top \\
			B & 0
		\end{bmatrix}^{-1}
		\approx
		\begin{bmatrix}
			I & -\tilde{A}^{-1}B^\top \\
			0 & I
		\end{bmatrix}
		\begin{bmatrix}
			\tilde{A}^{-1} & 0\\
			0 & \tilde{S}^{-1}
		\end{bmatrix}
		\begin{bmatrix}
			I & 0 \\
			-B\tilde{A}^{-1} & I
		\end{bmatrix},
	\end{equation}
	where $\tilde{A}^{-1}$ and $\tilde{S}^{-1}$ are approximations of the inverses of $A$ and of the Schur complement $S := -BA^{-1}B^\top$, respectively; note that equation \eqref{eq:block_factors} becomes exact by choosing $\tilde{A}^{-1} = A^{-1}$ and $\tilde{S}^{-1}=S^{-1}$. As noted in the works \cite{BO.2006,FMW.2019} in the context of the flow of a Navier-Stokes fluid, the (usually challenging) approximation of the Schur complement becomes extremely simple if one considers the modified system:
		\begin{equation}\label{eq:AL}
		\begin{bmatrix}
			A + \gamma B^\top M_p^{-1}B & B^\top \\
			B & 0
		\end{bmatrix}
		\begin{bmatrix}
\bz_h \\ p_h
		\end{bmatrix}
		=
		\begin{bmatrix}
			\bm{F} \\
			g
		\end{bmatrix},
	\end{equation}
	where $\gamma>0$ is called the augmented Lagrangian parameter and $M_p$ is the mass matrix of the pressure space $Q^h$. The new term $\gamma B^\top M_p^{-1}B \bz_h$ corresponds to a term $\gamma \nabla\diver \bu_h$, which could be interpreted as arising from an $L^2$-penalisation of the divergence. This modification does not change the solution (since the velocity $\bu_h$ is divergence-free), but it does change the Schur complement, allowing the very simple approximation $\tilde{S}^{-1} = - \gamma M_p^{-1}$ (see \cite{FMW.2019} for details). Note that for our choice of pressure space, $M_p$ is diagonal and can be inverted exactly.

	The main difficulty becomes the approximation of the inverse of the new stress-velocity operator, since a new term with a large kernel has been introduced (the system degenerates for large $\gamma$). In particular, the standard multigrid method breaks down for such systems. Building on the works \cite{Sch.1999,Sch.1999b,LWXZ.2007}, a Reynolds-robust multigrid preconditioner was developed in \cite{BO.2006,FMW.2019} by employing multigrid smoothers and transfer operators that take into account the divergence constraint; in particular, the smoothers are of additive Schwarz-type based on patch solves with appropriately chosen patches. In the context of non-Newtonian implicitly constituted fluids, these ideas were then applied in \cite{FG.2020} for a discretisation based on the Scott--Vogelius element; while possessing several advantages, such as $H^1$-conformity and exact enforcement of the divergence constraint, the method developed in  \cite{FG.2020} can become computationally expensive due to large patch solves for the multigrid smoothing and requires non-standard transfer operators. In contrast, an advantage of the $H(\diver)$-based discretisation employed in this work is that additive Schwarz smoothers based on standard star patches and the usual transfer operators are appropriate; see \cite{LFM.2022} for more details. In summary, in the present work we apply the block factorisation preconditioning strategy from \cite{FG.2020}, but in tandem with the multigrid method used in \cite{LFM.2022} for the stress-velocity block.
 
\section{Flow past an airfoil}

\begin{figure}
\centering
	\includegraphics[width=0.7\textwidth]{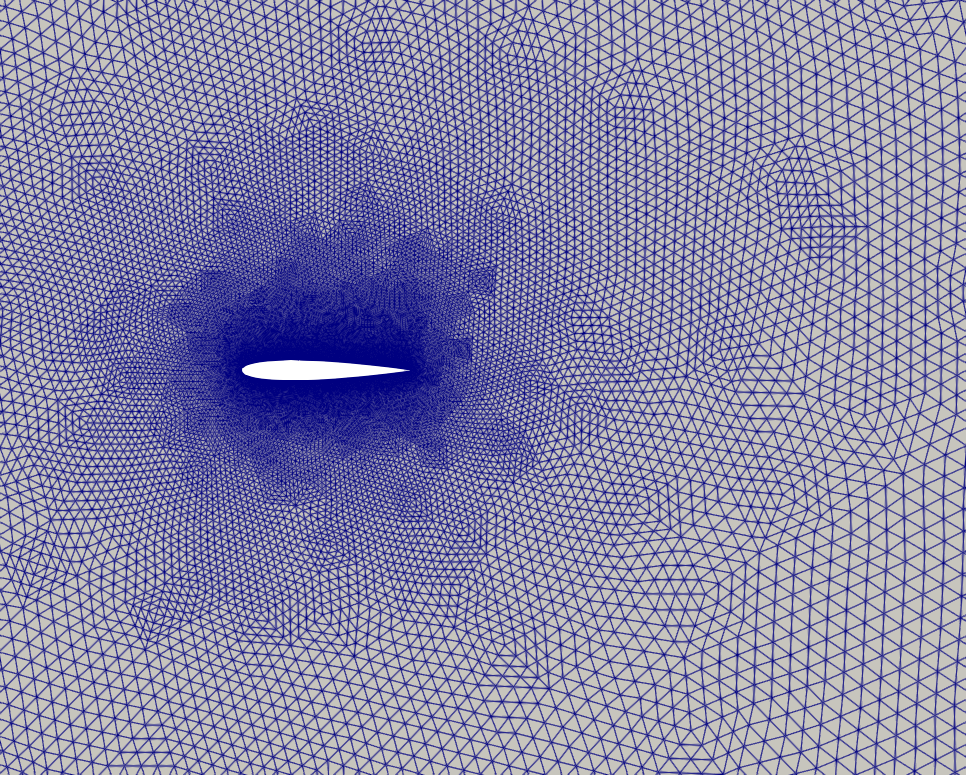}%
	\caption{Base mesh for the flow around an airfoil near the obstacle.}
	\label{fig:mesh}
\end{figure}

As a numerical example we will consider a square domain $(-8,8)^2$ with an airfoil of length 1 and maximum width 0.12 placed at the origin. The base mesh contains 8808 elements and is more refined towards the obstacle; a section of the mesh is shown in Figure \ref{fig:mesh}; the problem will then be solved on a uniform refinement of the base mesh with either $2.1\times 10^5$ total degrees of freedom (1 refinement) or $8.5\times 10^5$ degrees of freedom (2 refinements). On the left, top and bottom boundaries, a far-field velocity equal to $\bv_b = (10,0)^\top$ was prescribed; on the right boundary a natural boundary condition was imposed, which prescribes a pressure drop of magnitude 5.

For this problem we set $\Re = 500$ and employ either the activated Euler model with $\Eu = 15$ and $\varepsilon = 0.001$, or the Navier-Stokes model ($\Eu = 0$). We compare the solution for the flow past the airfoil for an activated Euler fluid against the solution to the full Navier-Stokes equations rather than the boundary layer equations and find very good agreement between the two solutions.
The example was implemented in \texttt{firedrake} \cite{Firedrake}.
The Newton linearisation was supplemented with the NLEQERR linesearch from PETSc \cite{PETSc}.



\begin{figure}
\centering
\subfloat[C][{\centering {\normalsize Magnitude of the velocity $\bv_h$}.}]{{%
	\includegraphics[width=0.9\textwidth]{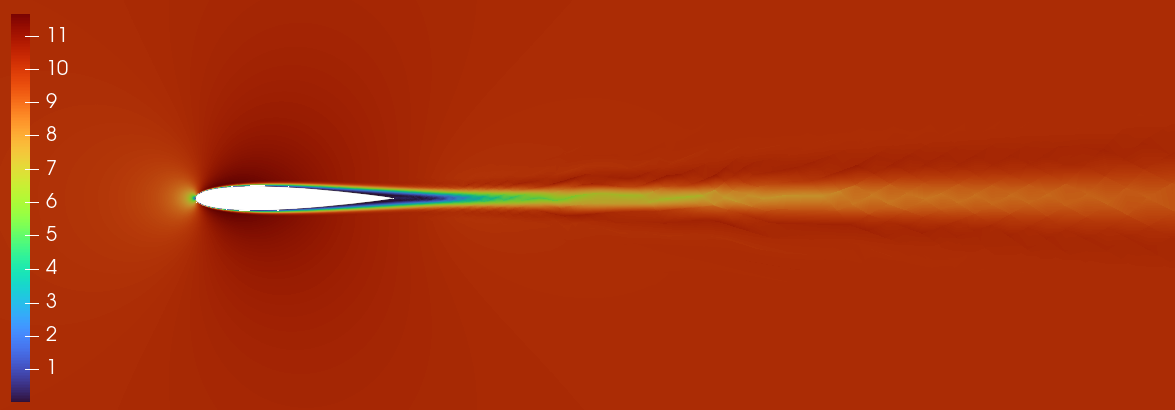}%
	}}\\
	\subfloat[D][{\centering {\normalsize Magnitude of the vorticity $|\omega_h|$}.}]{{%
	\includegraphics[width=0.9\textwidth]{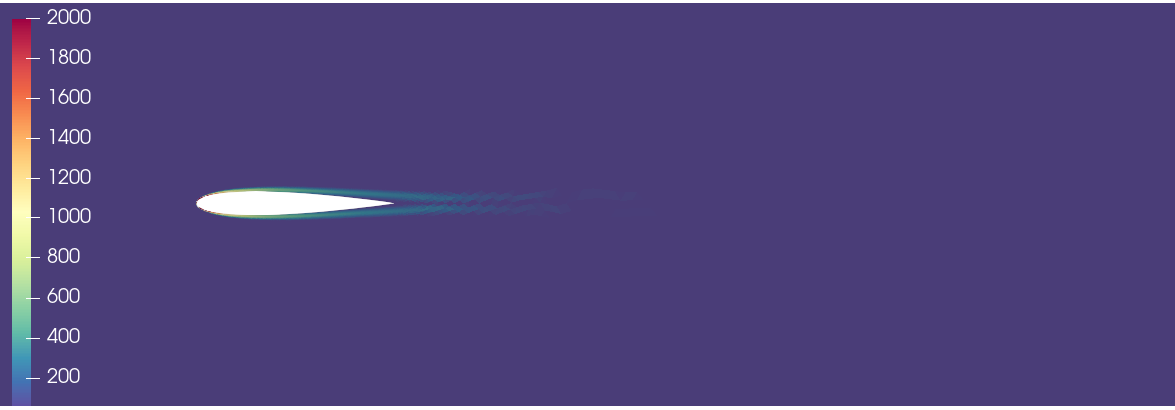}%
	}}\\
	\subfloat[E][{\centering {\normalsize Pressure $p_h$}.}]{{%
	\includegraphics[width=0.9\textwidth]{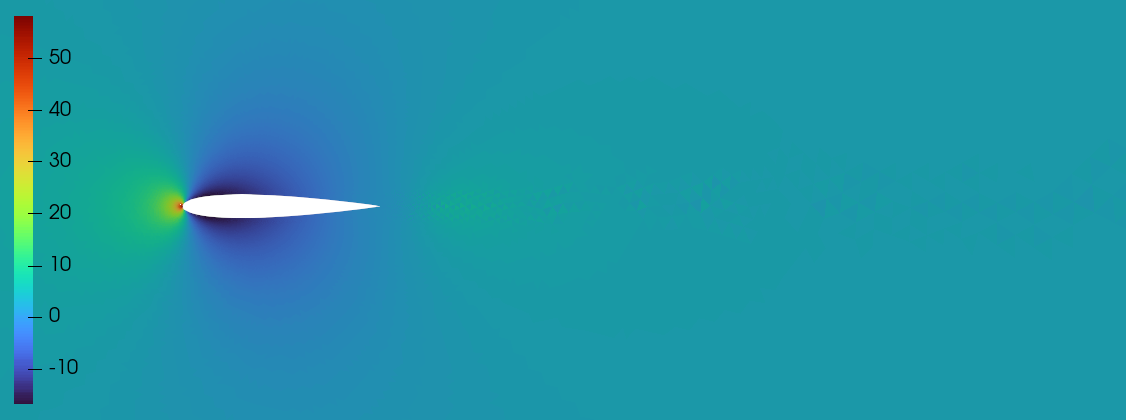}%
	}}\\
	\caption{Computed solution for the activated Euler model.}%
	\label{fig:solution}
\end{figure}

\begin{figure}
\centering
\subfloat[C][{\centering {\normalsize Magnitude of the velocity $\bv_h$}.}]{{%
	\includegraphics[width=0.9\textwidth]{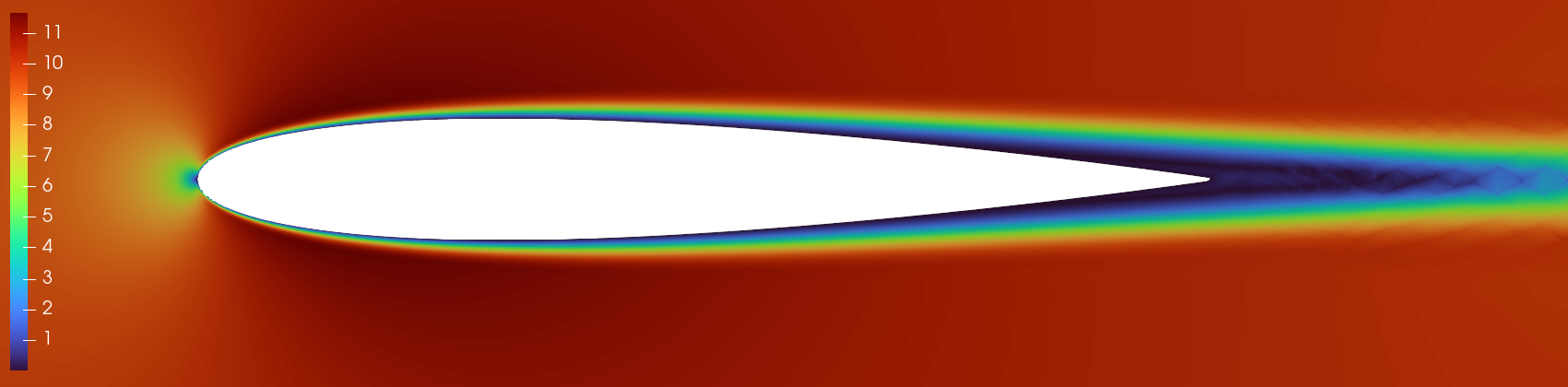}%
	}}\\
	\subfloat[D][{\centering {\normalsize Magnitude of the vorticity $|\omega_h|$}.}]{{%
	\includegraphics[width=0.9\textwidth]{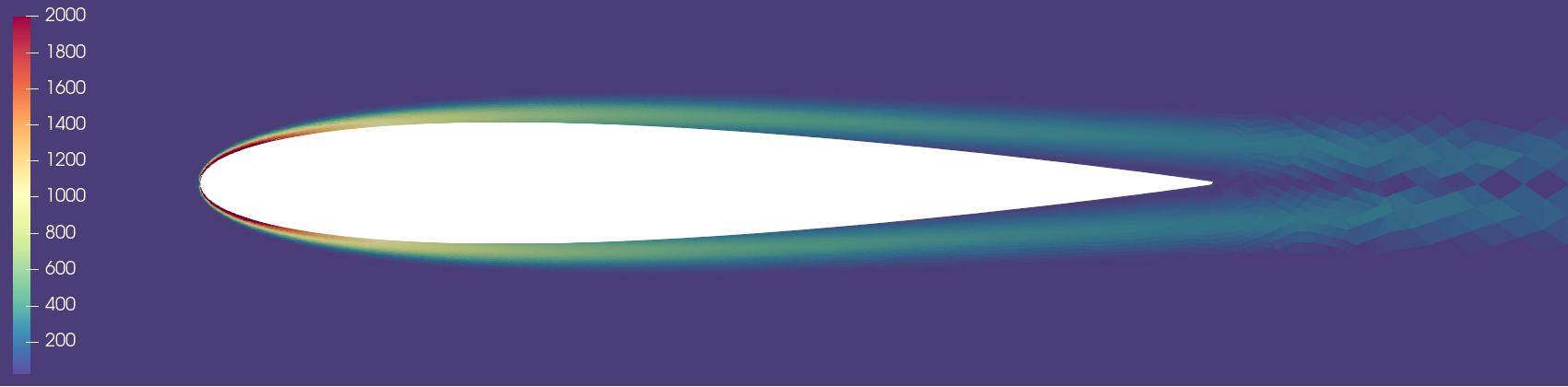}%
	}}\\
	\subfloat[E][{\centering {\normalsize Pressure $p_h$}.}]{{%
	\includegraphics[width=0.9\textwidth]{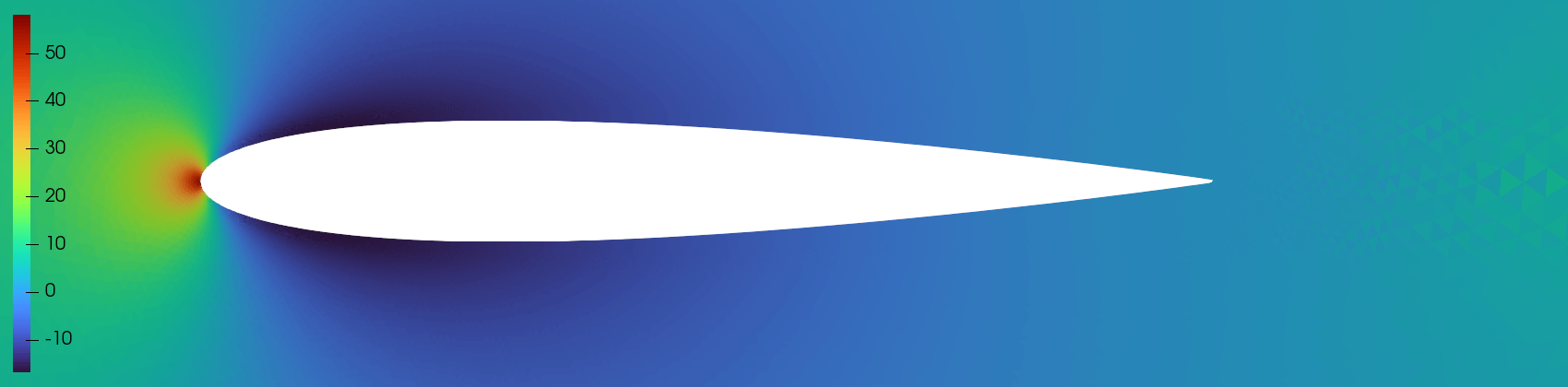}%
	}}\\
	\caption{Computed solution for the activated Euler model near the obstacle.}%
	\label{fig:solution_zoom}
\end{figure}

\begin{figure}
\centering
\subfloat[C][{\centering {\normalsize Full size}.}]{{%
	\includegraphics[width=0.5\textwidth]{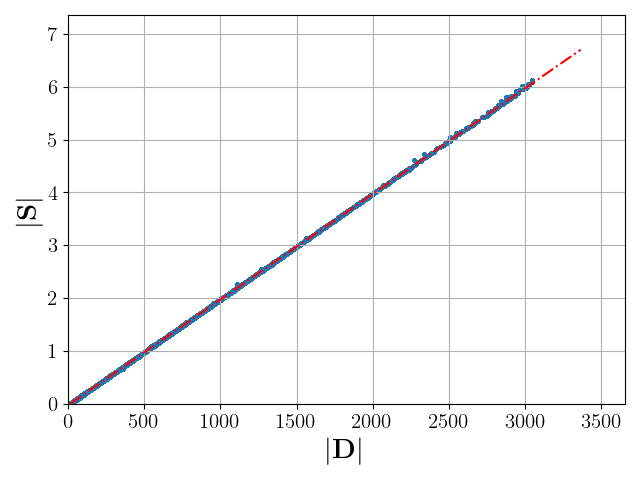}%
	}}%
	\subfloat[D][{\centering {\normalsize Zoom near the activation region}.}]{{%
	\includegraphics[width=0.5\textwidth]{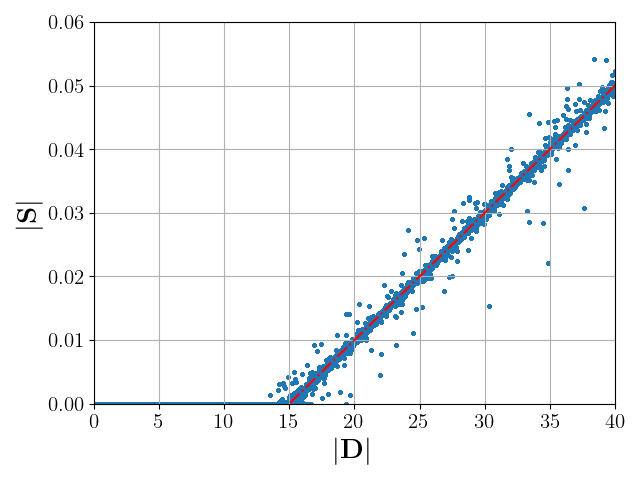}%
	}}\\
	\caption{Computed constitutive relation for the activated Euler model.}%
	\label{fig:CR_plot}
\end{figure}

\begin{figure}
\centering
	\includegraphics[width=0.7\textwidth]{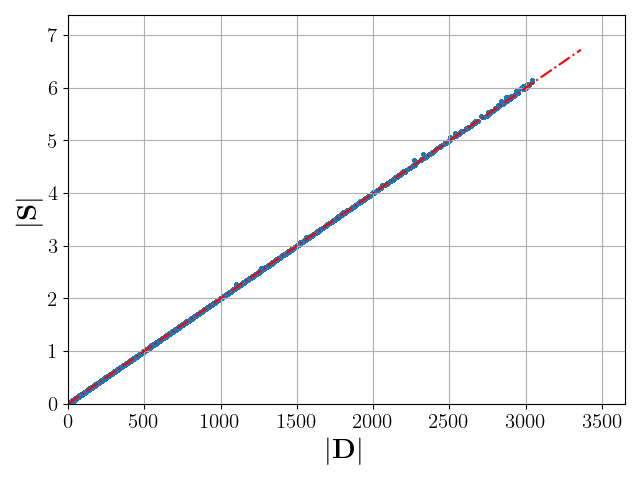}%
	\caption{Computed constitutive relation for the Navier-Stokes model.}
	\label{fig:CR_newtonian_plot}
\end{figure}

Figure \ref{fig:solution} shows the plots of the magnitude of the velocity, the magnitude of the vorticity $|\omega|:= |\partial_1 \bv_{h,2} - \partial_2 \bv_{h,1}|$, and the pressure on a portion of the domain near the obstacle;
Figure \ref{fig:solution_zoom} shows a close-up of the same solution (recall that the obstacle has length 1).
A plot of the relationship between $|\BS|$ and $|\BD|$ is shown in Figure \ref{fig:CR_plot} for the computed solutions corresponding to the activated Euler model on the second mesh refinement. The computed response relation between the norm of the stress and the norm of the symmetric part of the velocity gradient behaviour matches well with the activated Euler constitutive relation and acts as a check to our computations. The response for the Navier-Stokes constitutive relation is shown in Figure \ref{fig:CR_newtonian_plot}. It can be observed that the two response relations are strikingly similar.

\begin{figure}
\centering
	\subfloat[A][{\centering {\normalsize Velocity magnitude ($\Eu=0$)}.}]{{%
	\includegraphics[width=0.5\textwidth]{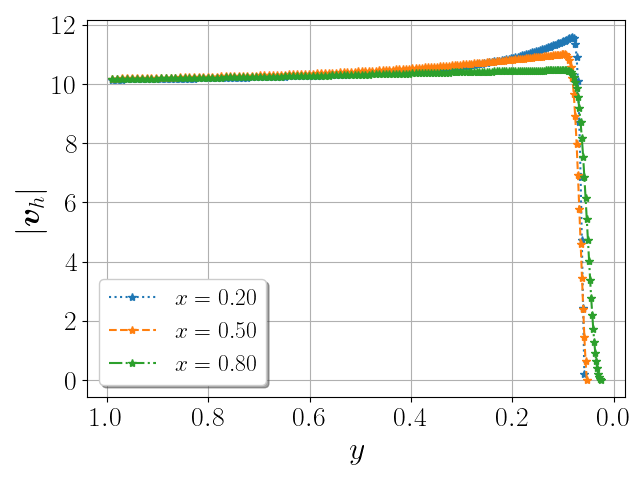}%
	}}%
	\subfloat[B][{\centering {\normalsize Velocity magnitude ($\Eu=15$)}.}]{{%
	\includegraphics[width=0.5\textwidth]{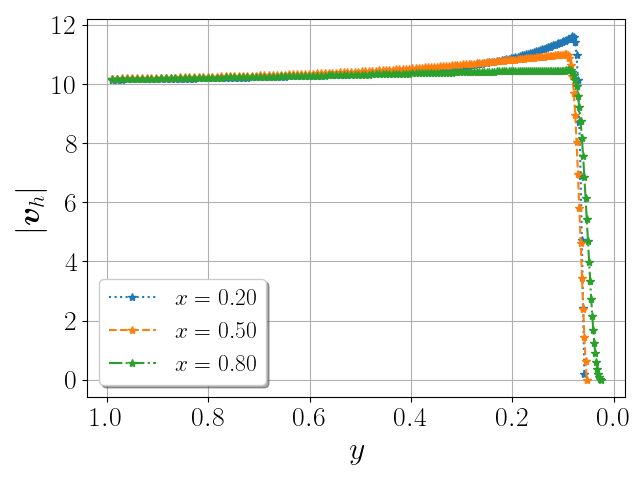}%
	}}\\
	\subfloat[A][{\centering {\normalsize Vorticity magnitude ($\Eu=0$)}.}]{{%
	\includegraphics[width=0.5\textwidth]{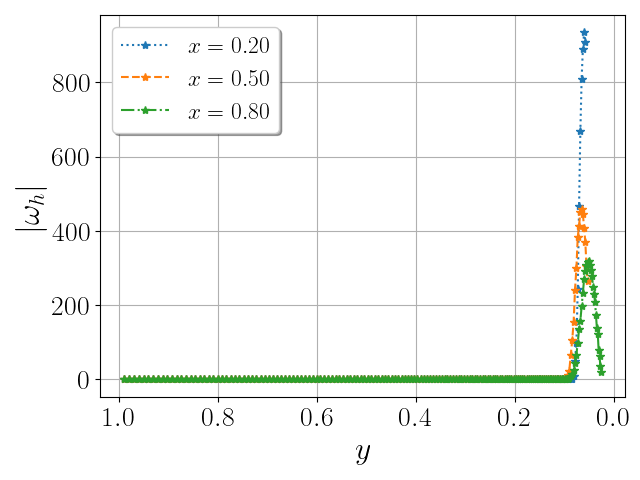}%
	}}%
	\subfloat[B][{\centering {\normalsize Vorticity magnitude ($\Eu=15$)}.}]{{%
	\includegraphics[width=0.5\textwidth]{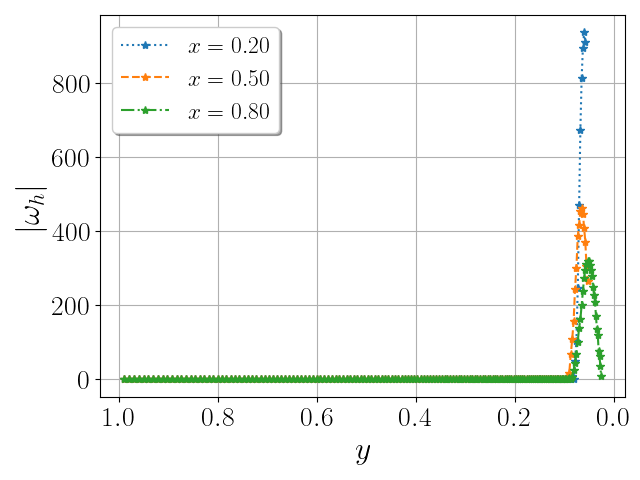}%
	}}\\
	\caption{Solutions corresponding to the Navier-Stokes (left) and activated Euler (right) models.}%
	\label{fig:blayer}
\end{figure}

In order to show the development of the boundary layer, we present in Figure \ref{fig:blayer} plots of the magnitude of the velocity $\bv_h$ and of the magnitude of the computed vorticity $\omega_h:=\frac{\partial}{\partial x}((\bv_h)_2) - \frac{\partial}{\partial y}((\bv_h)_1) $ along three vertical slices (choosing $x\in\{0.2,0.5,0.8\}$). Here again a striking resemblance between the Navier-Stokes (Newtonian) and the activated Euler (non-Newtonian) models can be observed.

\section{Concluding remarks}

This study, though short, has been able to establish several interesting results which we list below:
\begin{itemize}
    \item The computational analysis that has been carried out has clearly shown that the activated Euler fluid model has the same capability to capture the development of boundary layers as the Navier-Stokes fluid model when the flow takes place past a bluff body. The agreement of the results for the velocity and vorticity for the two fluids are strikingly similar.
    \item It is interesting that for the development of an efficient numerical scheme we found it convenient to view the studied model as a subclass of implicit constitutive relations, and also as the dual model to the one for a Bingham fluid. As the mathematical community has already dedicated significant attention to the analysis of the Bingham-type problems and to the construction of their efficient approximation, we could apply those tools to the activated Euler fluid. 
    \item With respect to the classical Euler equations, the mathematical properties exhibited by the governing equations for the activated Euler constitutive relation are remarkably better. For the model considered here, large data existence of weak solution has been established in \cite{BMR.2020} for both steady and unsteady flows and for different types of boundary conditions.
    \item We consider an implicit constitutive theory wherein the fluid exhibits distinct response characteristics below and above an activation criterion, that is, there is a bifurcation in the behavior of the fluid at a critical activation value. This is reminiscent of the response exhibited by rigid-plastic or elasto-plastic response of solids based on a yield criterion or solid-to-solid phase transition such as from Martensite to Austenite transition based on the criterion for maximal rate of entropy production.
 
    \item The activated Euler fluid model contains the activation parameter $\tau_*$ that, in this study, has to be specified a priori. There is however an interesting study  \cite{GP.2017B} where the identification of the activation criterion is part of the solution and it is determined in an aposteriori manner. 
\end{itemize}

\bibliographystyle{plain}
\bibliography{./bibliography}

\end{document}